# BOUNDEDNESS PROPERTIES OF FRACTIONAL INTEGRAL OPERATORS ASSOCIATED TO NON–DOUBLING MEASURES

JOSÉ GARCÍA-CUERVA AND A. EDUARDO GATTO

ABSTRACT. The main purpose of this paper is to investigate the behaviour of fractional integral operators associated to a measure on a metric space satisfying just a mild growth condition, namely that the measure of each ball is controlled by a fixed power of its radius. This allows, in particular, non–doubling measures. It turns out that this condition is enough to build up a theory that contains the classical results based upon the Lebesgue measure on euclidean space and their known extensions for doubling measures. We start by analyzing the images of the Lebesgue spaces associated to the measure. The Lipschitz spaces, defined in terms of the metric, play a basic role too. For a euclidean space equipped with one of these measures, we also consider the so-called "regular"*BMO* space introduced by X. Tolsa. We show that it contains the image of a Lebesgue space in the appropriate limit case and also that the image of the space "regular"*BMO* is contained in the adequate Lipschitz space.

## 1. INTRODUCTION

One of the most important tools of Harmonic Analysis in the last thirty years has been the notion of a space of homogeneous type. These spaces were formally introduced in [CW1], although in some form or another they were already present in some previous work (see, for example [CG]). The basic ingredients are a metric or quasi-metric space and a Borel measure $\mu$ satisfying the so–called "doubling condition", which means that there exists a constant $C$, such that, for every ball

*Date*: October 25, 2018.

1991 *Mathematics Subject Classification.* 2000 revision, 42B20, 26A33, 47B38, 47G10.

Supported in part by DGES, Spain, under grant BFM2001-0189 and, for the second author, by Sabbatical grants from the Program of the Spanish Ministry of Education and from DePaul University .

$B(x, r)$ of center $x$ and radius $r$

(1.1) $$\mu(B(x, 2r)) \leq C\mu(B(x, r)).$$

The success of the spaces of homogeneous type as the natural setting for a large portion of Harmonic Analysis, mainly the Calderón-Zygmund theory, led to the firm belief by almost all specialists in the field, that we had achieved the right level of generality to do Analysis. A measure of this success is the fact that the setting chosen by E. M. Stein in his book [S] to develop the basic theory is, essentially, that of a space of homogeneous type.

It came as a surprise to many when it was announced in [NTV2] that *"The doubling condition is superfluous for most of the classical theory"*. What they meant is that a rather complete theory of Calderón-Zygmund operators could be developed if, in a separable metric space, the condition (1.1) was replaced by the condition (2.1), which we introduce in the next section and which is generally referred to by saying that the measure is $n-$dimensional. A few recent articles pertaining to this line of research are [NTV1, NTV2, GM1, GM2, To1, To2, To3]

The present work is devoted to investigate the behavior of the fractional integral $I_\alpha$ associated to an $n-$dimensional measure $\mu$ on a metric space, as given by definition 3.1, for $0 < \alpha < n$ and other related operators $K_\alpha$ defined in section 4.

In section 3 we see that, for $1 < p < n/\alpha$, $I_\alpha$ maps $L^p(\mu)$ continuously into $L^q(\mu)$, where $1/q = 1/p - \alpha/n$, with a substitute weak–type result for $p = 1$. This generalizes the classical Hardy-Littlewood-Sobolev theorem. We also observe that condition (2.1) is actually necessary for this Hardy-Littlewood-Sobolev theorem to hold. These results, under slightly different conditions on the metric space, were also obtained by V. Kokilashvili and A. Meskhi in [KM].

In section 4 we introduce a general class of operators $K_\alpha$, which contain as particular cases the fractional integrals $I_\alpha$.

In section 5 we extend our interest to the case $p > n/\alpha$. Here the classical Lipschitz space defined in terms of the metric, which we shall denote $\mathcal{L}ip(\beta)$, plays a basic role, as one can see from the statements of the two theorems presented in section 5: Theorem 5.2, which says that when $p > n/\alpha$, $K_\alpha$ maps $L^p(\mu)$ continuously into $\mathcal{L}ip(\alpha - n/p)$; and Theorem 5.3, which establishes that $K_\alpha$ maps $\mathcal{L}ip(\beta)$ continuously into $\mathcal{L}ip(\alpha + \beta)$, provided it sends constants to constants.

Let us remark that, up to this point, we do not need to assume that the metric space is separable.

Finally, in section 6, we specialize the setting to some euclidean space $\mathbb{R}^d$ equipped with an $n$-dimensional measure $\mu$. We study the action of $K_\alpha$ on $L^{n/\alpha}(\mu)$ and show that this space is taken boundedly into $RBMO(\mu)$, the "regular"$BMO$ space introduced by X. Tolsa in [To3]. We also show in this section that the space $RBMO(\mu)$ is taken by $K_\alpha$ into $\mathcal{L}ip\,(\alpha)$, provided the right cancellation condition $K_\alpha(1) = 0$ holds.

The study of singular integral operators with respect to $\mu$ on the Lipschitz spaces will appear in a subsequent paper together with the characterization of these spaces in terms of mean oscillation.

## 2. The general setting. Basic facts

In this paper, unless otherwise stated, $(\mathbb{X}, d, \mu)$ will always be a metric measure space (that is, $d$ is a distance on $\mathbb{X}$ and $\mu$ is a Borel measure on $\mathbb{X}$), such that, for every ball

$$B(x,r) = \{y \in \mathbb{X} \,:\, d(x,y) < r\}, \quad x \in \mathbb{X}, \ r > 0,$$

we have

(2.1) $$\mu(B(x,r)) \leq C r^n,$$

where $n$ is some fixed positive real number and $C$ in independent of $x$ and $r$. Sometimes we shall refer to condition (2.1) by saying that the measure $\mu$ is $n-$dimensional.

Whenever we refer to "the ball $B$"in the sequel, we shall understand that we have chosen for it a fixed center and a fixed radius. That way, it makes sense to say that if $B$ is a ball and $k$ is a positive real number, we shall denote by $kB$ the ball having the same center as $B$ and radius $k$ times that of $B$.

For the sake of completeness we shall state and prove two very simple lemmas, which will be used all throughout the paper.

**Lemma 2.1.** *For every $\gamma > 0$*

(2.2) $$\int_{B(x,r)} \frac{1}{d(x,y)^{n-\gamma}} \,\mathrm{d}\mu(y) \leq C r^\gamma$$

**Proof.** If $n \leq \gamma$, (2.2) follows immediately from (2.1).

If $\gamma < n$, we write

$$\int_{B(x,r)} \frac{1}{d(x,y)^{n-\gamma}} \, d\mu(y) = \sum_{j=0}^{\infty} \int_{2^{-j-1}r \leq d(x,y) < 2^{-j}r} \frac{1}{d(x,y)^{n-\gamma}} \, d\mu(y)$$

$$\leq \sum_{j=0}^{\infty} \frac{1}{(2^{-j-1}r)^{n-\gamma}} \mu(B(x, 2^{-j}r)) \leq \sum_{j=0}^{\infty} \frac{2^{(j+1)(n-\gamma)}}{r^{n-\gamma}} C(2^{-j}r)^n$$

$$= C \sum_{j=0}^{\infty} 2^{-\gamma j} r^\gamma = Cr^\gamma.$$

$\square$

**Lemma 2.2.** *For every $\gamma > 0$*

(2.3) $$\int_{\mathbb{X} \setminus B(x,r)} \frac{1}{d(x,y)^{n+\gamma}} \, d\mu(y) \leq Cr^{-\gamma}$$

**Proof.**

$$\int_{\mathbb{X} \setminus B(x,r)} \frac{1}{d(x,y)^{n+\gamma}} \, d\mu(y) = \sum_{j=0}^{\infty} \int_{2^j r \leq d(x,y) < 2^{j+1}r} \frac{1}{d(x,y)^{n+\gamma}} \, d\mu(y)$$

$$\leq \sum_{j=0}^{\infty} \frac{\mu(B(x, 2^{j+1}r))}{(2^j r)^{n+\gamma}} \leq C \sum_{j=0}^{\infty} \frac{(2^{j+1}r)^n}{(2^j r)^{n+\gamma}}$$

$$= C \sum_{j=0}^{\infty} 2^{-\gamma j} r^{-\gamma} = Cr^{-\gamma}.$$

$\square$

## 3. Fractional integrals and the Hardy-Littlewood-Sobolev theorem

**Definition 3.1.** Let $0 < \alpha < n$. The fractional integral $I_\alpha$ associated to the measure $\mu$ will be defined, for appropriate functions $f$ on $\mathbb{X}$ as

(3.1) $$I_\alpha f(x) = \int_{\mathbb{X}} \frac{f(y)}{d(x,y)^{n-\alpha}} \, d\mu(y).$$

According to Lemma 2.1, for fixed $x$, the function $y \mapsto \dfrac{1}{d(x,y)^{n-\alpha}}$ is locally integrable with respect to $\mu$. Therefore, definition (3.1) makes perfectly good sense when $f$ is bounded and has bounded support. However, many times we shall need to define the fractional integral for larger classes of functions. In such cases, we shall explain in detail in the next sections how to carry out the corresponding extensions.

**Theorem 3.2.** *For $1 \leq p < \frac{n}{\alpha}$ and $\frac{1}{q} = \frac{1}{p} - \frac{\alpha}{n}$, we have*

$$(3.2) \qquad \mu\left(\{x \in \mathbb{X} \,:\, |I_\alpha f(x)| > \lambda\}\right) \leq \left(\frac{C \|f\|_{L^p(\mu)}}{\lambda}\right)^q,$$

*that is, $I_\alpha$ is a bounded operator from $L^p(\mu)$ into the Lorentz space $L^{q,\infty}(\mu)$.*

**Proof.** We are going to adapt to our context the proof given by Stein [S] for $\mathbb{R}^n$. We can take $f \geq 0$.

$$I_\alpha f(x) = \int_\mathbb{X} \frac{f(y)}{d(x,y)^{n-\alpha}} \, d\mu(y) = I + II,$$

where $I$ is the integral over $B(x, r)$ and $II$ is the integral over $\mathbb{X} \setminus B(x, r)$. By Hölder's inequality, if $p > 1$ and $p'$ denotes its conjugate exponent (i. e. $1/p + 1/p' = 1$).

$$|II| \leq \|f\|_{L^p(\mu)} \left(\int_{\mathbb{X} \setminus B(x,r)} \frac{1}{d(x,y)^{(n-\alpha)p'}} \, d\mu(y)\right)^{1/p'}.$$

$(n - \alpha)p' = n + \gamma$, where $\gamma = n(p' - 1) - \alpha p'$, so that

$$\frac{\gamma}{p'} = n\left(1 - \frac{1}{p'}\right) - \alpha = \frac{n}{p} - \alpha > 0.$$

By Lemma 2.2 $|II| \leq \|f\|_{L^p(\mu)} \left(Cr^{-\gamma}\right)^{1/p'} = C \|f\|_{L^p(\mu)} r^{-\left(\frac{n}{p} - \alpha\right)}$,

which holds even for $p = 1$. We can and do assume that $\|f\|_{L^p(\mu)} = 1$. Also, for $\lambda > 0$, we choose $r$ so that
$Cr^{-\left(\frac{n}{p} - \alpha\right)} = \lambda/2$. Then

$$\{x \in \mathbb{X} \,:\, |I_\alpha f(x)| > \lambda\} \subset \{x \in \mathbb{X} \,:\, |I| > \lambda/2\} \cup \{x \in \mathbb{X} \,:\, |II| > \lambda/2\}.$$

By the relation between $r$ and $\lambda$, the second of these sets is empty. We use Hölder's inequality once more to obtain

$$\begin{aligned}
|I| &\leq \left(\int_{B(x,r)} \frac{|f(y)|^p}{d(x,y)^{n-\alpha}} \, d\mu(y)\right)^{1/p} \left(\int_{B(x,r)} \frac{d\mu(y)}{d(x,y)^{n-\alpha}}\right)^{1/p'} \\
&\leq Cr^{\alpha/p'} \left(\int_{B(x,r)} \frac{|f(y)|^p}{d(x,y)^{n-\alpha}} \, d\mu(y)\right)^{1/p},
\end{aligned}$$

where we have also used Lemma 2.1. Then, by applying Tchebichev's inequality, we get

$$\mu(\{x \in \mathbb{X} : |I_\alpha f(x)| > \lambda\}) \leq \mu(\{x \in \mathbb{X} : |I| > \lambda/2\})$$
$$\leq Cr^{\alpha p/p'}\lambda^{-p} \int_{\mathbb{X}} \int_{B(x,r)} \frac{|f(y)|^p}{d(x,y)^{n-\alpha}} \, d\mu(y) \, d\mu(x)$$
$$= Cr^{\alpha p/p'}\lambda^{-p} \int_{\mathbb{X}} \int_{B(y,r)} \frac{d\mu(x)}{d(x,y)^{n-\alpha}} |f(y)|^p \, d\mu(y)$$
$$\leq Cr^{\alpha p/p'} r^\alpha \lambda^{-p} = Cr^n = C\lambda^{-q}, \text{ since } \lambda = Cr^{-(n/p-\alpha)}.$$

□

**Corollary 3.3.** *For $1 < p < \frac{n}{\alpha}$ and $\frac{1}{q} = \frac{1}{p} - \frac{\alpha}{n}$, we have*

(3.3) $$\|I_\alpha f\|_{L^q(\mu)} \leq C \|f\|_{L^p(\mu)}$$

**Proof.** It suffices to apply Marcinkiewicz's interpolation theorem with indices slightly bigger and slightly smaller than $p$.

□

The next result reveals that condition (2.1) is the minimal requirement one must have in order for (3.3) or (3.2) to be valid. Note that we consider only measures without atoms, so that the fractional integral given by (3.1) is well defined.

**Theorem 3.4.** *For a measure $\mu$, finite over balls and not having any atoms, condition (2.1) is necessary for the Hardy-Littlewood-Sobolev theorem to hold.*

**Proof.** Suppose that (3.3) holds. Let $B$ be a ball of radius $r$. If $\mu(B) = 0$, then (2.1) is trivially true. Let $\mu(B) \neq 0$. For each $x \in B$, we have

$$I_\alpha \chi_B(x) = \int_B \frac{1}{d(x,y)^{n-\alpha}} \, d\mu(y) \geq \frac{1}{(2r)^{n-\alpha}} \mu(B).$$

By applying (3.3) we get

$$\frac{1}{(2r)^{n-\alpha}} \mu(B)^{1+\frac{1}{q}} \leq \left(\int_B |I_\alpha \chi_B(x)|^q \, d\mu(x)\right)^{1/q}$$
$$\leq C \|\chi_B\|_{L^p(\mu)} = C\mu(B)^{1/p},$$

which is equivalent to

(3.4) $$\mu(B)^{1+\frac{1}{q}-\frac{1}{p}} \leq C (r^n)^{1-\frac{\alpha}{n}}.$$

Since $1 + \frac{1}{q} - \frac{1}{p} = 1 - \frac{\alpha}{n}$, inequality (3.4) is, precisely, condition (2.1). A similar argument works if we assume (3.2). □

## 4. Fractional integral operators

Now we shall associate to our fixed $n-$dimensional measure, a family of fractional kernels and the corresponding fractional integral operators.

**Definition 4.1.** Let $0 < \alpha < n$ and $0 < \varepsilon \leq 1$. A function $k_\alpha : \mathbb{X} \times \mathbb{X} \longrightarrow \mathbb{C}$ is said to be a fractional kernel of order $\alpha$ and regularity $\varepsilon$ if it satisfies the following two conditions:

$$(4.1) \qquad |k_\alpha(x,y)| \leq \frac{C}{d(x,y)^{n-\alpha}}, \quad \text{for all } x \neq y;$$

and

$$(4.2) \qquad |k_\alpha(x,y) - k_\alpha(x',y)| \leq C\frac{d(x,x')^\varepsilon}{d(x,y)^{n-\alpha+\varepsilon}}.$$

for $d(x,y) \geq 2\,d(x,x')$. The corresponding operator $K_\alpha$, which will be called "fractional integral operator", will be given by

$$(4.3) \qquad K_\alpha(f)(x) = \int_\mathbb{X} k_\alpha(x,y)\,f(y)\,d\mu(y).$$

By (4.1), $K_\alpha(f)$ is well defined for $f \in L^p(\mu)$, $1 \leq p < \frac{n}{\alpha}$ and (3.3) or (3.2) are also valid for it, as for $I_\alpha(f)$. Next we see that the fractional integral $I_\alpha$ is an example of fractional integral operador with a kernel having regularity 1.

**Lemma 4.2.** Let $x, y, z \in \mathbb{X}$ be such that $2d(x,y) \leq d(x,z)$. Then

$$(4.4) \qquad \left| \frac{1}{d(x,z)^{n-\alpha}} - \frac{1}{d(y,z)^{n-\alpha}} \right| \leq C\frac{d(x,y)}{d(x,z)^{n-\alpha+1}}.$$

**Proof.** By the mean value theorem of real differential calculus, we have, for $s, t > 0$

$$\left| s^{\alpha-n} - t^{\alpha-n} \right| \leq (n-\alpha)\left|(1-\theta)s + \theta t\right|^{\alpha-n-1} |s-t|$$

for some $\theta \in\,]0,1[$. Now, since $2d(x,y) \leq d(x,z)$, we get

$$\left| \frac{1}{d(x,z)^{n-\alpha}} - \frac{1}{d(y,z)^{n-\alpha}} \right| \leq \frac{C\,|d(x,z) - d(y,z)|}{d(x,z)^{n-\alpha+1}} \leq \frac{Cd(x,y)}{d(x,z)^{n-\alpha+1}}.$$

$\square$

**Definition 4.3.** Let $k_\alpha$ be a fractional kernel of order $\alpha$ and regularity $\varepsilon$, $f \in L^p(\mu)$, $p > n/\alpha$, and $\alpha - \frac{n}{p} < \varepsilon$. We define

$$(4.5) \qquad \widetilde{K_\alpha}f(x) = \int_\mathbb{X} \{k_\alpha(x,y) - k_\alpha(x_0,y)\}\,f(y)\,\mathrm{d}\mu(y),$$

where $x_0$ is some fixed point of $\mathbb{X}$.

We observe that the integral in (4.5) converges both locally and at $\infty$ as a consequence of (4.1), (4.2) and Hölder's inequality. Of course the function just defined depends on the election of $x_0$. But the difference between any two functions obtained in (4.5) for different elections of $x_0$ is just a constant.

## 5. Lipschitz spaces

From now on, we shall assume that $\mu(\mathbb{X}) = \infty$. The results below are true also when $\mu(\mathbb{X}) < \infty$, but other results are more appropriate in that case. They will be treated elsewhere.

**Definition 5.1.** Given $\beta \in ]0,1[$, we shall say that the function $f : \mathbb{X} \to \mathbb{C}$ satisfies a Lipschitz condition of order $\beta$ provided

$$(5.1) \qquad |f(x) - f(y)| \leq C d(x,y)^\beta \text{ for every } x, y \in \mathbb{X}$$

and the smallest constant in inequality (5.1) will be denoted by $\|f\|_{\mathcal{L}ip(\beta)}$ It is easy to see that the linear space of all Lipschitz functions of order $\beta$, modulo constants, becomes, with the norm $\|\ \|_{\mathcal{L}ip(\beta)}$, a Banach space, which we shall call $\mathcal{L}ip(\beta)$.

**Theorem 5.2.** *Let $k_\alpha$ be a fractional kernel with regularity $\varepsilon$. If $n/\alpha < p \leq \infty$ and $\alpha - \dfrac{n}{p} < \varepsilon$, then $\widetilde{K_\alpha}$ maps $L^p(\mu)$ boundedly into $\mathcal{L}ip\left(\alpha - \dfrac{n}{p}\right)$.*

**Proof.** Assume, first $p < \infty$. Consider $x \neq y$ and let $B$ be the ball with center $x$ and radius $r = d(x,y)$. Then, we have

$$\left|\widetilde{K_\alpha}f(x) - \widetilde{K_\alpha}f(y)\right| \leq \int_{2B} |k_\alpha(x,z)| |f(z)| \, d\mu(z)$$
$$+ \int_{2B} |k_\alpha(y,z)| |f(z)| \, d\mu(z) + \int_{\mathbb{X}\setminus 2B} |k_\alpha(x,z) - k_\alpha(y,z)| |f(z)| \, d\mu(z).$$

We shall estimate each of these three terms separately. For the first two terms we use (4.1) and Hölder's inequality.

$$\int_{2B} |k_\alpha(x,z)| |f(z)| \, d\mu(z)$$
$$\leq \int_{2B} \frac{|f(z)|}{d(x,z)^{n-\alpha}} \, d\mu(z) \leq \|f\|_{L^p(\mu)} \left(\int_{2B} \frac{d\mu(z)}{d(x,z)^{(n-\alpha)p'}}\right)^{1/p'}$$

Now observe that $(n-\alpha)p' = n - p'\left(\alpha - \dfrac{n}{p}\right)$ and, since $\alpha - \dfrac{n}{p} > 0$, the integral converges and, by Lemma 2.1, we have

$$\int_{2B} |k_\alpha(x,z)|\,|f(z)|\,\mathrm{d}\mu(z) \leq C\,\|f\|_{L^p(\mu)}\,(2r)^{\alpha-(n/p)}.$$

The second term is estimated in a similar way after noting that $2B \subset B(y, 3r)$.

Next, to estimate the third term we use (4.2) and Hölder's inequality to get

$$\int_{\mathbb{X}\setminus 2B} |k_\alpha(x,z) - k_\alpha(y,z)|\,|f(z)|\,\mathrm{d}\mu(z)$$

$$\leq \int_{\mathbb{X}\setminus 2B} \frac{C d(x,y)^\varepsilon}{d(x,z)^{n-\alpha+\varepsilon}}\,|f(z)|\,\mathrm{d}\mu(z)$$

$$\leq C d(x,y)^\varepsilon \|f\|_{L^p(\mu)} \left(\int_{\mathbb{X}\setminus 2B} \frac{\mathrm{d}\mu(z)}{d(x,z)^{(n-\alpha+\varepsilon)p'}}\right)^{1/p'}.$$

Note that $(n - \alpha + \varepsilon)p' = n + p'\left(\dfrac{n}{p} + \varepsilon - \alpha\right)$ and since, by hypothesis $\dfrac{n}{p} + \varepsilon - \alpha > 0$, the integral converges and, by Lemma 2.2

$$\int_{\mathbb{X}\setminus 2B} |k_\alpha(x,z) - k_\alpha(y,z)|\,|f(z)|\,\mathrm{d}\mu(z) \leq C\,\|f\|_{L^p(\mu)}\,d(x,y)^{\alpha-\frac{n}{p}}.$$

Putting together the three estimates, we get

$$\left|\widetilde{K_\alpha}f(x) - \widetilde{K_\alpha}f(y)\right| \leq C\,\|f\|_{L^p(\mu)}\,d(x,y)^{\alpha-\frac{n}{p}},$$

as we wanted to prove.

The case $p = \infty$ is similar, only easier. All the estimates proved for $p < \infty$ continue to hold for $p = \infty$.

$\square$

We shall show next that the Lipschitz spaces with respect to the metric are preserved by the fractional integral operators provided that the image of a constant is also constant.

**Theorem 5.3.** *Let $k_\alpha$ be a fractional kernel with regularity $\varepsilon$, and $\alpha, \beta > 0$ be such that $\alpha + \beta < \varepsilon$. Then $\widetilde{K_\alpha}$ is a bounded mapping from $\mathcal{L}ip\,(\beta)$ to $\mathcal{L}ip\,(\alpha + \beta)$ if and only if $\widetilde{K_\alpha}(1) = 0$.*

**Proof.** To see that the condition is necessary, note that continuity of the operator $\widetilde{K_\alpha}$ implies that $\widetilde{K_\alpha}(1)$ must be constant. On the other hand $\widetilde{K_\alpha}(1)(x_0) = 0$. Therefore, the constant has to be 0.

To prove the sufficiency we consider $x \neq y$ points of $\mathbb{X}$, and we want to estimate $\left|\widetilde{K_\alpha}(f)(x) - \widetilde{K_\alpha}(f)(y)\right|$.

Observe first that

$$\widetilde{K_\alpha}(1) = 0 \Leftrightarrow \widetilde{K_\alpha}(1)(x) - \widetilde{K_\alpha}(1)(y) = 0$$
$$\Leftrightarrow \int_\mathbb{X} \{k_\alpha(x,z) - k_\alpha(y,z)\} \, \mathrm{d}\mu(z) = 0,$$

where the integrals above converge because $0 < \alpha < \varepsilon$.

Thus we can write

$$\widetilde{K_\alpha}(f)(x) - \widetilde{K_\alpha}(f)(y)$$
$$= \int_\mathbb{X} \{k_\alpha(x,z) - k_\alpha(y,z)\} (f(z) - f(x)) \, \mathrm{d}\mu(z) = I + II,$$

where $I$ is the integral over $2B$, $B$ being the ball with center $x$ and radius $r = d(x,y)$ and $II$ is the integral over $\mathbb{X} \setminus 2B$. Now

$$|I| \leq \int_{2B} \frac{1}{d(x,z)^{n-\alpha}} |f(z) - f(x)| \, \mathrm{d}\mu(z) + \int_{2B} \frac{1}{d(y,z)^{n-\alpha}} |f(z) - f(x)| \, \mathrm{d}\mu(z)$$

In the sum above, both terms can be estimated in the same fashion. For example, for the first one, using Lemma 2.1 we get

$$\int_{2B} \frac{d(x,z)^\beta \, \mathrm{d}\mu(z)}{d(x,z)^{n-\alpha}} \leq C(2r)^{\alpha+\beta} \leq Cd(x,y)^{\alpha+\beta}$$

and for the second one, extending the integral to $B(y, 3r)$, we get the same estimate.

In order to estimate $II$, we use (4.2) and Lemma 2.2 to obtain

$$|II| \leq C \int_{\mathbb{X} \setminus 2B} \frac{d(x,y)^\varepsilon \, d(x,z)^\beta}{d(x,z)^{n-\alpha+\varepsilon}} \, \mathrm{d}\mu(z) \leq Cd(x,y)^\varepsilon \int_{\mathbb{X} \setminus 2B} \frac{\mathrm{d}\mu(z)}{d(x,z)^{n+\varepsilon-\alpha-\beta}}$$
$$\leq Cd(x,y)^\varepsilon r^{\alpha+\beta-\varepsilon} \leq Cd(x,y)^{\alpha+\beta}.$$

This finishes the proof. □

## 6. The "regular" *BMO* space of X. Tolsa

For this section, the underlying metric measure space will be $\mathbb{X} = \mathbb{R}^d$ with an $n-$dimensional measure $\mu$, $n \leq d$. It is for this context for which X. Tolsa has introduced (in [To3]) the space "regular"*BMO*, although his definition makes sense, in principle, in our general setting.

All balls considered in this section will be assumed to be centered at points of the support of $\mu$.

**Definition 6.1.** A function $f \in L^1_{\text{loc}}(\mu)$ has "regular bounded mean oscillation" with respect to $\mu$ if the following two conditions are satisfied, where $\rho > 1$ is a fixed constant.

  a) There exists a constant $C$, such that, for every ball $B$

$$\text{(6.1)} \qquad \int_B |f(x) - m_B(f)| \, d\mu(x) \leq C\mu(\rho B),$$

where $m_B(f) = \dfrac{1}{\mu(B)} \displaystyle\int_B f(x) \, d\mu(x)$ ; and also

  b) there is a constant $C$, such that, for every two balls $B \subset V$, $B$ having radius $r$ and $V$ having radius $s$ :

$$\text{(6.2)} \qquad |m_B(f) - m_V(f)| \leq C K_{B,V} \left( \frac{\mu(\rho B)}{\mu(B)} + \frac{\mu(\rho V)}{\mu(V)} \right),$$

where

$$\text{(6.3)} \qquad K_{B,V} = 1 + \sum_{k=1}^{N_{B,V}} \frac{\mu(2^k B)}{(2^k r)^n},$$

$N_{B,V}$ being the first integer $k$ such that $2^k r \geq s$.

Denoting by $\|f\|_\star$ the smallest constant which makes both (6.1) and (6.2) true, the space $RBMO(\mu)$ obtained by considering equal those functions of regular bounded mean oscillation with respect to $\mu$ that differ by a constant, becomes a Banach space with the norm $\| \, \|_\star$. It is a basic fact, proved by X. Tolsa, that the space does not depend on the constant $\rho > 1$ used. Besides, Tolsa has shown that, in the definition of $RBMO(\mu)$ one can use cubes rather than balls, obtaining exactly the same space.

Before stating and proving the first result of this section, we need two lemmas. The first one is an extension to our setting of the one given in [GV2] for the doubling case.

**Lemma 6.2.** *Let $f \in L^{n/\alpha}(\mu)$ be a function vanishing outside $\rho B$, where $B$ is a ball and $\rho \geq 1$. Then*

$$\text{(6.4)} \qquad \int_B |I_\alpha f(x)| \, d\mu(x) \leq C \, \|f\|_{L^{n/\alpha}(\mu)} \, \mu(\rho B)$$

*with $C$ independent of $f$, $\rho$ and $B$.*

**Proof.** Choose $p \in ]1, n/\alpha[$ and let $1/q = 1/p - \alpha/n$. Then, by applying Jensen's inequality (twice) and the Hardy-Littlewood-Sobolev theorem 3.3, we get

$$\frac{1}{\mu(\rho B)} \int_B |I_\alpha f(x)| \, d\mu(x) \leq \left(\frac{1}{\mu(\rho B)} \int_{\rho B} |I_\alpha f(x)|^q \, d\mu(x)\right)^{1/q}$$

$$\leq \frac{C}{\mu(\rho B)^{1/q-1/p}} \left(\frac{1}{\mu(\rho B)} \int_{\rho B} |f(x)|^p \, d\mu(x)\right)^{1/p}$$

$$\leq \frac{C}{\mu(\rho B)^{1/q-1/p}} \left(\frac{1}{\mu(\rho B)} \int_{\rho B} |f(x)|^{n/\alpha} \, d\mu(x)\right)^{\alpha/n} = C \, \|f\|_{L^{n/\alpha}(\mu)}$$

□

**Lemma 6.3.** *let $B \subset V$ be balls. Then*
$$V \subset 2^{N_{B,V}+1} B \subset 5V,$$
*where $N_{B,V}$ is the integer appearing in (6.3).*

**Proof.** Denote the centers of $B$ and $V$ by $x_B$ and $x_V$ respectively and let the corresponding radii be $r$ and $s$. By the definition of $N_{B,V}$, we have $2^{N_{B,V}-1}r < s \leq 2^{N_{B,V}}r$. Then

$$y \in V \Rightarrow d(y, x_V) < s \Rightarrow d(y, x_B) \leq d(y, x_V) + d(x_V, x_B) < 2s \leq 2^{N_{B,V}+1}r.$$

So, we have seen that $V \subset 2^{N_{B,V}+1}B$. Also

$$y \in 2^{N_{B,V}+1}B \Rightarrow d(y, x_B) < 2^{N_{B,V}+1}r \Rightarrow$$
$$d(y, x_V) \leq d(y, x_B) + d(x_B, x_V) < 2^{N_{B,V}+1}r + s \leq 5s.$$

This proves that $2^{N_{B,V}+1}B \subset 5V$.

□

Next, as in the case of Lipschitz functions that we considered in last section, we have to redefine $K_\alpha$ on $L^{n/\alpha}$. Note that the elements of *RBMO* are classes of functions modulo constants.

**Definition 6.4.** Let $0 < \alpha < \varepsilon \leq 1$ and $k_\alpha(x,y)$ be a fractional kernel of order $\alpha$ with regularity $\varepsilon$. For $f \in L^{n/\alpha}(\mu)$, we define

$$(6.5) \quad \overline{K_\alpha}(f)(x) = \int_{\mathbb{X}} \left\{k_\alpha(x,y) - \chi_{\mathbb{X} \setminus B(x_0,1)} k_\alpha(x_0, y)\right\} f(y) \, d\mu(y)$$

for some fixed $x_0 \in \mathbb{X}$. We shall show below in Theorem 6.5 that $\overline{K_\alpha} f(x)$ is well defined almost everywhere with respect to $\mu$.

Although the definition depends on the choice of $x_0$, ; different choices of $x_0$ yield functions that differ only by a constant.

**Theorem 6.5.** *Let $0 < \alpha < \varepsilon \leq 1$ and $f \in L^{n/\alpha}(\mu)$. Then $\overline{K_\alpha}(f)$ is well defined almost everywhere with respect to $\mu$, $\overline{K_\alpha}(f) \in RBMO(\mu)$ and*
$$\left\|\overline{K_\alpha}(f)\right\|_\star \leq C \left\|f\right\|_{L^{n/\alpha}(\mu)}$$
*with $C$ independent of $f$.*

**Proof.** We shall show first condition (6.1) and, at the same time, the existence almost everywhere of the integral in (6.5). It suffices to show that, for each ball $B = B(x_1, r)$, there is a constant $c_B$ such that

$$(6.6) \qquad \int_B \left|\overline{K_\alpha}(f)(x) - c_B\right| \, d\mu \leq C \|f\|_{L^{n/\alpha}(\mu)} \mu(2B).$$

Let
$$c_B = \int_{\mathbb{X}} \left\{\chi_{\mathbb{X} \setminus B(x_1, 2r)}(y) k_\alpha(x_1, y) - \chi_{\mathbb{X} \setminus B(x_0, 1)}(y) k_\alpha(x_0, y)\right\} f(y) \, d\mu(y).$$

Note that $|K_\alpha(f)(x) - c_B|$ is dominated by:

$$\int_{\mathbb{X}} \left|k_\alpha(x, y) - \chi_{\mathbb{X} \setminus B(x_1, 2r)}(y) k_\alpha(x_1, y)\right| |f(y)| \, d\mu(y) = I(x) + II(x),$$

where $I(x)$ is the integral over $2B$ and $II(x)$ the integral over the complement of $2B$. Next we estimate the integrals over $B$ of $I(x)$ and $II(x)$.

$$\int_B I(x) \, d\mu(x) \leq \int_B \int_{\mathbb{X}} |k_\alpha(x,y)| \chi_{2B}(y) |f(y)| \, d\mu(y) \, d\mu(x)$$
$$\leq \int_B |I_\alpha(\chi_{2B}|f|)| \, d\mu(x) \leq C\|f\|_{L^{n/\alpha}(\mu)} \mu(2B),$$

where the last inequality follows from Lemma 6.2. To estimate the integral of $II(x)$ over $B$ observe first that, since $x \in B$ and $y \in \mathbb{X} \setminus 2B$ using (4.2),

$$II(x) \leq \int_{\mathbb{X} \setminus 2B} \frac{d(x, x_1)^\varepsilon}{d(x_1, y)^{n-\alpha+\varepsilon}} |f(y)| \, d\mu(y)$$

and, by Hölder's inequality and Lemma 2.2 we get that $II(x) \leq C\|f\|_{L^{n/\alpha}(\mu)}$. Therefore the integral of $II(x)$ over $B$ is also bounded by $C\|f\|_{L^{n/\alpha}(\mu)} \mu(2B)$.

Next we shall prove (6.2). Let $B \subset V$ be balls, $B$ having radius $r$ and $V$ radius $s$. We shall show

$$(6.7) \quad \frac{1}{\mu(B)} \frac{1}{\mu(V)} \int_B \int_V \left|\overline{K_\alpha} f(x) - \overline{K_\alpha} f(y)\right| \, d\mu(x) \, d\mu(y)$$
$$\leq C K_{B,V} \|f\|_{L^{n/\alpha}(\mu)} \left(\frac{\mu(\rho B)}{\mu(B)} + \frac{\mu(\rho V)}{\mu(V)}\right)$$

with $\rho = 10$, say.

Observe that the left hand side of (6.7) dominates
$$\left|m_B(\overline{K_\alpha}f) - m_V(\overline{K_\alpha}f)\right|,$$
so that we get (6.2) for $\overline{K_\alpha}(f)$.

Let $2^{N_{B,V}+2}B = \widetilde{V}$ for which we know from lemma 6.3, that $2V \subset \widetilde{V} \subset 10V$.

For almost every $x$ and $y$, we can write
$$\overline{K_\alpha}(f)(x) - \overline{K_\alpha}(f)(y) = \int_{\mathbb{X}} (k_\alpha(x,z) - k_\alpha(y,z))\, f(z)\, \mathrm{d}\mu(z)$$

Then
$$\left|\overline{K_\alpha}(f)(x) - \overline{K_\alpha}(f)(y)\right| \leq I_\alpha(|f|\chi_{2B})(x) + I_\alpha(|f|\chi_{\widetilde{V}\setminus 2B})(x) + I_\alpha(|f|\chi_{\widetilde{V}})(y)$$
$$+ \int_{\mathbb{X}\setminus\widetilde{V}} |k_\alpha(x,z) - k_\alpha(y,z)|\,|f(z)|\,\mathrm{d}\mu(z).$$

Let us denote by $I, II, III$ and $IV$ the four terms in the right hand side of the last inequality. We estimate separately the double average of each of these four terms. The first and the third are dealt with by means of the Lemma 6.2.

$$\frac{1}{\mu(B)}\frac{1}{\mu(V)}\int_B\int_V I\,\mathrm{d}\mu(y)\,\mathrm{d}\mu(x)$$
$$= \frac{1}{\mu(B)}\int_B I_\alpha(|f|\chi_{2B})(x)\,\mathrm{d}\mu(x) \leq C\,\|f\|_{n/\alpha}\,\frac{\mu(2B)}{\mu(B)}$$

and, likewise

$$\frac{1}{\mu(B)}\frac{1}{\mu(V)}\int_B\int_V III\,\mathrm{d}\mu(y)\,\mathrm{d}\mu(x)$$
$$= \frac{1}{\mu(V)}\int_V I_\alpha(|f|\chi_{\widetilde{V}})(y)\,\mathrm{d}\mu(y) \leq C\,\|f\|_{n/\alpha}\,\frac{\mu(10V)}{\mu(V)}.$$

In order to deal with $II$, we use Hölder's inequality, obtaining

$$II = \int_{\widetilde{V}\setminus 2B} \frac{1}{d(x,z)^{n-\alpha}} |f(z)|\,\mathrm{d}\mu(z)$$
$$\leq \left(\int_{\widetilde{V}\setminus 2B} \frac{\mathrm{d}\mu(z)}{d(x,z)^n}\right)^{(n-\alpha)/n} \|f\|_{n/\alpha} \leq CK_{B,V}\,\|f\|_{n/\alpha},$$

since

$$\int_{\widetilde{V}\setminus 2B} \frac{\mathrm{d}\mu(z)}{d(x,z)^n} \leq \sum_{k=1}^{N_{B,V}+1} \int_{2^{k+1}B\setminus 2^k B} \frac{\mathrm{d}\mu(z)}{d(x,z)^n} \leq CK_{B,V}.$$

and, finally, by (4.2)

$$IV \leq Cd(x,y)^\varepsilon \int_{\mathbb{X}\setminus 2V} \frac{1}{d(x,z)^{n-\alpha+\varepsilon}} |f(z)| \, \mathrm{d}\mu(z)$$

$$\leq Cd(x,y)^\varepsilon \|f\|_{L^{n/\alpha}(\mu)} \left( \int_{\mathbb{X}\setminus 2V} \frac{\mathrm{d}\mu(z)}{d(x,z)^{(n-\alpha+\varepsilon)n/(n-\alpha)}} \right)^{(n-\alpha)/n}.$$

Since $(n-\alpha+\varepsilon)\dfrac{n}{n-\alpha} = n + \dfrac{\varepsilon n}{n-\alpha}$, applying Lemma 2.2, we get

$$IV \leq Cs^\varepsilon \|f\|_{L^{n/\alpha}(\mu)} \left( s^{-\varepsilon n/(n-\alpha)} \right)^{(n-\alpha)/n} = C \|f\|_{L^{n/\alpha}(\mu)}.$$

This finishes the proof.

□

Finally we present a result which can be viewed either as an extension of the case $p = \infty$ of Theorem 5.2 or as an extension of the case $\beta = 0$ of Theorem 5.3.

**Theorem 6.6.** *Let $k_\alpha$ be a fractional kernel of order $\alpha$ and regularity $\epsilon$, with $0 < \alpha < \epsilon$. Then $\widetilde{K_\alpha}(f)$, as given by (4.5), makes sense for any function in $RBMO(\mu)$ and the operator $\widetilde{K_\alpha}$ thus defined, is bounded from $RBMO(\mu)$ to $\mathcal{L}ip\,(\alpha)$ if and only if $\widetilde{K_\alpha}(1) = 0$.*

**Proof.** The necessity follows exactly as in the proof of Theorem 5.3. Namely, the boundedness of $\widetilde{K_\alpha}(f)$ implies that $\widetilde{K_\alpha}(1)$ must be a constant. But, since $\widetilde{K_\alpha}(1)(x_0) = 0$ by definition, it follows that $\widetilde{K_\alpha}(1) = 0$ or, equivalently, that

$$(6.8) \qquad \int_{\mathbb{X}} (k_\alpha(x,y) - k_\alpha(x_0,y)) \, \mathrm{d}\mu(y) = 0$$

for all $x$.

To prove the sufficiency, we will use the following equivalent definition of $RBMO(\mu)$ given in [To3]. $f \in RBMO(\mu)$ if and only if, there are constants $C$, and also a constant $f_U$ for every ball $U$, such that the following two conditions hold

$$(6.9) \qquad \frac{1}{\mu(\rho U)} \int_U |f(x) - f_U| \, \mathrm{d}\mu(x) \leq C$$

with some fixed $\rho$ and also, for each couple of balls $U \subset W$

$$(6.10) \qquad |f_U - f_W| \leq CK_{U,W}.$$

In addition, the infimum of the constants $C$ in (6.9) and (6.10) is equivalent to $\|f\|_\star$.

Let $f \in RBMO(\mu)$. Take two points $x \neq y$ and let $B = B(x,r)$ with $r = d(x,y)$. Then

$$\widetilde{K_\alpha}(f)(x) - \widetilde{K_\alpha}(f)(y) = \int_{\mathbb{X}} (k_\alpha(x,z) - k_\alpha(y,z)) f(z) \, d\mu(z)$$

$$= \int_{\mathbb{X}} (k_\alpha(x,z) - k_\alpha(y,z)) (f(z) - f_{2B}) \, d\mu(z),$$

where we have used the cancellation property (6.8).

We can bound this difference by breaking up the absolute value of the integrand in the following way

$$\left|\widetilde{K_\alpha}(f)(x) - \widetilde{K_\alpha}(f)(y)\right| \leq I_\alpha \left(|f - f_{2B}| \chi_{2B}\right)(x) + I_\alpha \left(|f - f_{2B}| \chi_{2B}\right)(y)$$

$$+ \int_{\mathbb{R}^d \setminus 2B} |k_\alpha(x,z) - k_\alpha(y,z)| \, |f(z) - f_{2B}| \, d\mu(z) = I + II + III.$$

We estimate each of these three terms separately.

The first term is estimated exactly as we did with the first term in the proof of Theorem 5.2, namely

$$I = \int_{2B} \frac{1}{d(x,z)^{n-\alpha}} |f(z) - f_{2B}| \, d\mu(z)$$

$$\leq \left(\int_{2B} \frac{d\mu(z)}{d(x,z)^{(n-\alpha)p'}}\right)^{1/p'} \left(\int_{2B} |f(z) - f_{2B}|^p \, d\mu(z)\right)^{1/p}$$

$$\leq C r^{\alpha - \frac{n}{p}} \mu(\rho B)^{\frac{1}{p}} \|f\|_\star \leq C r^\alpha \|f\|_\star,$$

where we have used Hölder's inequality with some $p > \frac{n}{\alpha}$, Lemma 2.1 and the fact that the space $RBMO(\mu)$ can also be defined with $L^p$ norms for $p > 1$ instead of just $L^1$ norms, a fact that follows from the version of the John-Nirenberg theorem established by Tolsa for $RBMO(\mu)$ (see [To3]).

Actually, the second term can be dealt with exactly in the same way as the first one, after noting that $2B \subset B(y, 3r)$. We obtain, just as before, $II \leq Cr^\alpha \|f\|_\star$.

So, all that is left is to examine the third term. For that we use the regularity of the kernel, since $d(x,y) = r < \dfrac{d(x,z)}{2}$. We have, taking into account that the constants $K_{B,V}$ can be bounded by an absolute constant when $B$ and $V$ have radii which are "comparable"(i. e. with ratio bounded above and below by two fixed positive numbers)

$$III \leq \int_{\mathbb{R}^d \setminus 2B} \frac{d(x,y)^\epsilon}{d(x,z)^{n-\alpha+\epsilon}} |f(z) - f_{2B}| \, \mathrm{d}\mu(z) \leq C d(x,y)^\epsilon$$

$$\times \sum_{k=1}^\infty \left\{ \int_{2^{k+1}B \setminus 2^k B} \frac{|f(z) - f_{2^{k+1}B}|}{d(x,z)^{n-\alpha+\epsilon}} \, \mathrm{d}\mu(z) + \frac{|f_{2^{k+1}B} - f_{2B}|}{(2^k r)^{n-\alpha+\epsilon}} \mu(2^{k+1}B) \right\}$$

$$\leq C r^\epsilon \|f\|_\star \left\{ \sum_{k=1}^\infty \frac{\mu(\rho 2^{k+1}B)}{(2^k r)^{n-\alpha+\epsilon}} + \sum_{k=1}^\infty k \frac{\mu(2^{k+1}B)}{(2^k r)^{n-\alpha+\epsilon}} \right\}$$

$$= C r^\epsilon \left\{ \sum_{k=1}^\infty (2^k r)^{\alpha-\epsilon} + \sum_{k=1}^\infty k (2^k r)^{\alpha-\epsilon} \right\} \|f\|_\star = C \|f\|_\star r^\alpha.$$

$\square$

DPTO. DE MATEMÁTICAS, C-XV, UNIVERSIDAD AUTÓNOMA, 28049, MADRID, SPAIN

*E-mail address*: `jose.garcia-cuerva@uam.es`

DEPARTMENT OF MATHEMATICS, DEPAUL UNIVERSITY, CHICAGO, ILLINOIS, 60614, U. S. A.

*E-mail address*: `aegatto@condor.depaul.edu`